\documentclass[reqno,12pt]{amsart}

\usepackage{amssymb,amsmath,amsthm,color}
\usepackage{mathtools}

\newtheorem{theorem}{Theorem}

\newtheorem{lemma}{Lemma}

\makeatletter
\newcommand*{\rom}[1]{\expandafter\@slowromancap\romannumeral #1@}

\author{Jing-Jing Huang}
\address{
JJH: Department of Mathematics and Statistics, University of Nevada, Reno,
1664 N. Virginia St., Reno, NV 89557}
\email{jingjingh@unr.edu}
\author{Huixi Li}
\address{
HL: Department of Mathematics and Statistics, University of Nevada, Reno,
1664 N. Virginia St., Reno, NV 89557}
\email{huixil@unr.edu}

\dedicatory{}

\begin{document}

\title{On a Generalization of a Theorem of Popov}
\makeatletter
\def\@setauthors{%
  \begingroup
  \def\thanks{\protect\thanks@warning}%
  \trivlist
  \centering\footnotesize \@topsep30\p@\relax
  \advance\@topsep by -\baselineskip
  \item\relax
  \author@andify\authors
  \def\\{\protect\linebreak}%
  \authors%
  \ifx\@empty\contribs
  \else
    ,\penalty-3 \space \@setcontribs
    \@closetoccontribs
  \fi
  \endtrivlist
  \endgroup
}
\def\@settitle{\begin{center}%
  \baselineskip14\p@\relax
    \bfseries
  \@title
  \end{center}%
}
\makeatother

\maketitle

\begin{abstract}
In this paper, we obtain sharp estimates for the number of lattice points under and near the dilation of a general parabola, the former generalizing  an old result of Popov. We   apply Vaaler's lemma and the Erd\H{o}s-Turan inequality to reduce the two underlying counting problems to mean values of a certain quadratic exponential sums, whose treatment is subject to classical analytic techniques.
\end{abstract}

\section{Introduction}

The Gauss circle problem and the Dirichlet divisor problem are both long standing open problems in analytic number theory. 
They ask for the best possible estimates when counting the number of lattice points inside the dilated circle $x^2+y^2=a^2$ or under the dilated hyperbola $xy=a^2$ as $a$ approaches infinity. Since there are only three types of conics, that is ellipse, hyperbola, and parabola, one may as well ask the same question for the parabola. 

The first result was obtained by Popov \cite{Popov} in this regard.  Let $\mathcal{P}_0$ be the standard parabola $y=x^2$, then $y = \frac{x^2}{a}$ is the dilation of $\mathcal{P}_0$ under the transformation $(x,y)\to\left(\frac{x}a,\frac{y}a\right)$. Our goal is to count the number of lattice points $(m,n)\in\mathbb{Z}^+\times\mathbb{Z}^+$ 
with $m\le b$ that are on or under the curve  $y = \frac{x^2}{a}$, i.e. to estimate the sum
$$
\sum_{1\leq x \leq b} \left\lfloor \frac{x^2}{a} \right\rfloor.
$$
As is customary for many \emph{rounding error problems},   we expect that the sum $\sum\lfloor \theta_n\rfloor-\theta_n+1/2$ should be small when the sequence $\{\theta_n\}$ is reasonably well distributed $\bmod 1$.
As such, we  introduce the notation 
$$
 E_{\mathcal{P}_0}(a, b)=\sum_{1 \leq x \leq b} \left\lfloor \frac{ x^2}{a} \right\rfloor
-
\sum_{1 \leq x \leq b}
\left( \frac{ x^2}{a} - \frac{1}{2} \right)
$$
to represent the error term of the corresponding lattice points counting problem for the parabola $\mathcal{P}_0$.   

Popov's result states that for large real numbers $a$ and $b$ such that $b\ll a$, we have
\begin{equation}\label{e1.1}
E_{\mathcal{P}_0}(a, b)=O \left(
a^{\frac{1}{2}} \exp \left( 
\frac{\left(\frac{3}{4} \ln 2 + \varepsilon \right) \log a}{\log \log a}
\right)
\right). 
\end{equation}
Throughout the paper $\exp(x)$ means $e^x$. 
He also made a remark that the error term would be the same if we consider the parabola $y = x^2 + k x$ instead of $y = x^2$ for some positive constant $k$. 

As Popov points out, the exponent $1/2$ in the error term of \eqref{e1.1} is best possible in general and can be regarded as achieving the \emph{square root cancellation}. 
Note that for both the Gauss circle problem and the Dirichlet divisor problem,  error terms of the same quality are still wide open conjectures. We only mention that recently Bourgain and Watt \cite{BourgainWatt} obtained the best known error $O\left(a^\frac{517}{824}\right)$ for both problems.

Nevertheless, it transpires that the bound \eqref{e1.1} is susceptible to further improvement. In \cite{HuangLi}, the authors prove that for positive integers $a$ and $b$, we have
\begin{equation}\label{e1.3}
E_{\mathcal{P}_0}(a, b)=O \left( a^{\frac{1}{2}}\log a+b a^{-\frac{1}{2}}\exp\left(\frac{(2+\varepsilon)(\log a)^{\frac{1}{2}}}{ \log\log a} \right)\right).
\end{equation}
It is readily seen that when $a$ is an integer and $b\ll a$ the above bound is indeed better than that of Popov \eqref{e1.1}.

On the opposite end, the minimum order of the error is studied in \cite{Chamizo}, in which Chamizo and Pastor prove that 
 there exist infinitely many positive integers $a$ such that 
\[ |E_{\mathcal{P}_0}(a, a)|
\geq C a^{\frac{1}{2}} \exp\left(\frac{(\sqrt{2} - \varepsilon) (\log a)^{\frac{1}{2}}} {\log \log a}\right)\]
for some absolute positive constant $C$. 
This shows that \eqref{e1.3} is surprisingly close to the true order of $E_{\mathcal{P}_0}(a, b)$.

The first main result of this paper is a generalization of Popov's theorem to an arbitrary parabola $\mathcal{P}: y=\alpha x^2+\beta x+\gamma$. Note that the dilation of $\mathcal{P}$ under the transformation $(x,y)\to\left(\frac{x}a,\frac{y}a\right)$ is $a \mathcal{P}: y=\frac{\alpha x^2}{a}+\beta x+\gamma a$. 
\begin{theorem}\label{t1} Let 
$a>1$, $b > 1$, $\alpha\neq0$, $\beta$ and $\gamma$ be real numbers.  Then
\[
\sum_{1 \leq x \leq b} \left\lfloor \frac{\alpha x^2}{a}+\beta x+\gamma a \right\rfloor
=
\sum_{1\le x \leq b}
\left( \frac{\alpha x^2}{a}+\beta x+\gamma a - \frac{1}{2} \right)
+ 
E_\mathcal{P}(a, b), 
\]
where
\[
E_{\mathcal{P}} (a, b)
\ll 
\left(a^{\frac{1}{2}} + b a^{- \frac{1}{2}}\right) \exp\left( \frac{\left(\frac{1}{2}\ln 2 + \varepsilon \right) \log \Delta}{\log \log \Delta } \right)
\]
and $\Delta =\max\left( \frac{b^2 a^{\frac{1}{2}}}{a + b},3\right)$. 
\end{theorem}

When $\mathcal{P}$ lies above the $x$-axis (i.e. $\alpha>0$ and $\beta^2\le 4\alpha\gamma$), we may similarly interpret  the sum $\sum_{1 \leq x \leq b} \left\lfloor \frac{\alpha x^2}{a}+\beta x+\gamma a  \right\rfloor$ as the number of lattice points $(m, n) \in \mathbb{Z}^+\times\mathbb{Z}^+$ with $m\le b$ that are on or under the dilated parabola $a \mathcal{P}$, or equivalently, the number of points of the form $\left( \frac{m}{a}, \frac{n}{a}\right)$ that are on or under the parabola $\mathcal{P}$, where $(m, n) \in \mathbb{Z}^+\times\mathbb{Z}^+$, $ m \leq b$. 

Compared to \eqref{e1.1}, we do not require $b \ll a$ in Theorem \ref{t1}, which holds for all $b$. Actually, we see that $E_{\mathcal{P}} (a, b)\ll a^{\frac12+\varepsilon}$ when $b\le a$ and $E_{\mathcal{P}} (a, b)\ll \frac{b^{1+\varepsilon}}{\sqrt{a}}$ when $b\ge a$.  Moreover, it is easy to check that the error term in Theorem \ref{t1} is slightly better than that of Popov's result \eqref{e1.1} when $b \ll a^{1-\varepsilon}$ for some $\varepsilon>0$. 

A problem of similar flavor is to estimate the number of lattice points close to the parabola. For $\delta\in\left(0,\frac{1}{2}\right)$, let
\[
A_\mathcal{P}(a, b, \delta) 
= \sum_{\substack{1 \leq x \leq b \\ \left\| \frac{\alpha x^2}a+\beta x+\gamma a \right\| < \delta }} 1.
\]
Here and throughout the paper $\|x\|=\min_{n\in\mathbb{Z}}|x-n|$ denotes the distance from $x$ to a nearest integer. Geometrically, the function $A_\mathcal{P}(a, b, \delta)$ counts the number of lattice points $(m, n) \in \mathbb{Z}^+\times\mathbb{Z}^+$ with $m\le b$ that are close to $a\mathcal{P}$, or equivalently the number of points of the form $(\frac{m}a,\frac{n}a)$ that are close to $\mathcal{P}$, where $(m, n) \in \mathbb{Z}^+\times\mathbb{Z}^+$ and $m \leq b$.

In \cite{HuangLi} the authors have proved that for any positive integers $a$ and $b$ with $a \geq 3$ and $\delta\in(0,\frac{1}{2})$, we have
\begin{equation}\label{e1.2}
A_{\mathcal{P}_0}(a,b,\delta)=2\delta b+O\left( a^{\frac{1}{2}} \log a+b a^{-\frac{1}{2}}\exp\left(\frac{(2+\varepsilon)(\log a)^{\frac{1}{2}}}{ \log\log a} \right)\right).
\end{equation}

We are able to generalize this to an arbitrary parabola as well as to remove the requirement that $a$ and $b$ be integers.
\begin{theorem}\label{t2} Let $a>1$, $b > 1$, $\delta\in(0,1/2)$, $\alpha\neq0$, $\beta$ and $\gamma$ be real numbers. Then
\[
A_\mathcal{P}(a, b, \delta) 
=
2\delta b
+ E_\mathcal{P}(a, b, \delta), 
\]
where
\[
E_{\mathcal{P}} (a, b, \delta)
\ll 
\left(a^{\frac{1}{2}} + b a^{- \frac{1}{2}}\right) \exp\left( \frac{\left(\frac{1}{2}\ln 2 + \varepsilon \right) \log \Delta}{\log \log \Delta } \right)
\]
and $\Delta$ is the same as in Theorem \ref{t1}.
\end{theorem}
It is worth noting that in our previous work \cite{HuangLi} the proofs of \eqref{e1.3} and \eqref{e1.2} are based on an elegant elementary argument and in particular involve quadratic Gauss sums associated with modulus $a$. Therefore, the lack of integral condition on $a$ in the current paper is more than merely a technical nuisance. In fact,  our approach is completely different from that of \cite{HuangLi}. We begin with some Fourier analytic methods to transform the counting problem into one about mean values of a certain quadratic exponential sums, then we execute a \emph{Weyl differencing} argument to obtain essentially optimal bounds on the latter sums. 

We end our introduction by briefly mentioning that the general problem of estimating the number of integral/rational points close to a submanifold in the Euclidean space is a very active area of research  recently, see \cite{B, BDV, BVVZ1, BVVZ2, H1, H3, H4, H2, H5, HuangLiu, VV}. In particular, there are close connections between this area and metric diophantine approximation on manifolds. We recommend to the interested readers the monographs  \cite{BD, Sp} for the background and the survey \cite{BRV} for recent developments. 

Throughout the paper, we use the notation $e(x)=e^{2\pi i x}$. We use Vinogradov's symbol $f(x)\ll g(x)$ and Landau's symbol $f(x)=O(g(x))$ to mean there exists a constant $C$ such that $|f(x)|\le Cg(x)$. We use $\varepsilon$ to denote any sufficiently small, positive constant, which may not necessarily be the same in each occurrence. 

\section{The proof of Theorem \ref{t1}}
Let $\psi(x) =  \{x\} - 1/2$, 
where $\{x\}=x-\lfloor x\rfloor$ denotes the fractional part of $x$. With this notation,  we may write the error term in the form
\begin{align*}
E_\mathcal{P}(a, b)
= 
\left|\sum_{1\leq x\le b}\psi\left(\frac{\alpha x^2}{a} + \beta x + \gamma a \right)\right|.
\end{align*}
We need the following Lemma which is essentially due to Vaaler \cite[Theorem 18]{Vaaler}. 
\begin{lemma}[{\cite[Corollary 6.2\footnote{The summation there is over the range $(N, 2N]$, but it is straightforward to rewrite the lemma in the stated form by shifting the domain of $f$.}]{AT}}]\label{Lemma1}
Let $H$ be a positive integer and $f:[M, M+N]\to\mathbb{R}$ be any function. Then we have$$
\left|\sum_{M<n\le M+N}\psi(f(n))\right|\le \frac{N}{2H+2}+\left(1+\frac1\pi\right)\sum_{1\leq h \leq H}\frac1h\left|\sum_{M<n\le M+N}e(hf(n))\right|.
$$
\end{lemma}
We apply Lemma \ref{Lemma1} with $M=0$, $N=b$ and $f(x)=\frac{\alpha x^2}{a} + \beta x + \gamma a$.
It then follows that
\begin{align}
E_{\mathcal{P}}(a, b) 
\le&
\frac{b}{2H+2}+\left(1+\frac1\pi\right)\sum_{1\leq h \leq H}\frac1h\left|\sum_{1\leq x\le b}e\left( h \left( \frac{\alpha x^2}{a} + \beta x + \gamma a \right) \right)\right|.\label{e2.1}
\end{align}
Clearly, without loss of generality, we may assume $\alpha>0$, as the other case $\alpha<0$ amounts to taking complex conjugate of the sum inside the absolute value in \eqref{e2.1}.
By the Cauchy-Schwarz inequality, we obtain
\begin{align}
\nonumber&\sum_{1\leq h \leq H}
\frac1h
\left|\sum_{1 \leq x\le b}e\left(h\left(\frac{\alpha x^2}a+\beta x+\gamma a\right)\right)\right|\\
\le&\left(\sum_{1\leq h \leq H}
\frac1h\right)^{\frac{1}{2}}
\left(\sum_{1\leq h \leq H}\frac1h\left|\sum_{1 \leq x\le b}e\left(h\left(\frac{\alpha x^2}a+\beta x+\gamma a\right)\right)\right|^2\right)^{\frac{1}{2}}.\label{CS}
\end{align}
Note that
\begin{align}\label{secondsum}
&\sum_{1\leq h \leq H}\frac1h\left|\sum_{1 \leq x\le b}e\left(h\left(\frac{\alpha x^2}a+\beta x+\gamma a\right)\right)\right|^2 \nonumber \\
=&\sum_{1\leq h \leq H}\frac1h
\sum_{0\leq |l|<b}e\left(h\left(\frac{\alpha l^2}a+\beta l\right)\right)\sum_{1\le y, y+l\le b}e\left(\frac{2\alpha hl y}a\right) \nonumber \\
\le&\sum_{1\leq h \leq H}\frac1h
\sum_{0\leq |l|<b}
\left|\sum_{1\le y, y+l\le b}e\left(\frac{2\alpha h l y}a\right)\right| \nonumber \\
\le&\sum_{1\leq h \leq H}\frac1h
\left(b+\sum_{1\leq l<b}\min\left(2b,\left\|\frac{2\alpha h l}a\right\|^{-1}\right)\right) \nonumber \\
\ll &b\log (H)+\sum_{1\leq j<b H}\min\left(2b,\left\|\frac{2\alpha j}a\right\|^{-1}\right)\sum_{\substack{h|j\\h>\frac{j}{b}}}\frac1h \nonumber \\
\ll &b\log (H)+\sum_{1\leq j<b H}\min\left(2b,\left\|\frac{2\alpha j}a\right\|^{-1}\right)\sum_{\substack{h|j\\1 \leq h < b}}\frac{h}{j} \nonumber \\
\ll &b\log (H)+I, 
\end{align}
where
$$
I=\sum_{1\le j<b H}\min\left(2b,\left\|\frac{2\alpha j}a\right\|^{-1}\right)\frac{g(j)}j
$$
and
$$
g(j)=\sum_{\substack{h|j\\1\leq h<b}}h.
$$
Thus, by \eqref{e2.1}, \eqref{CS}, and \eqref{secondsum}, we have
\begin{equation}\label{e2.2}
E_{\mathcal{P}}(a, b) \ll \frac{b}{H} + (\log H)^{\frac{1}{2}} (b \log H + I)^{\frac{1}{2}} \ll \frac{b}{H} + b^{\frac{1}{2}} \log H + I^{\frac{1}{2}} (\log H)^{\frac{1}{2}}. 
\end{equation}

Now we estimate the sum $I$. Due to the behavior of the function $g(j)$, we write the interval $[1, b H)$ as the union of two sub-intervals $[1, b)$ and $[b, b H)$. By the bound on the divisor sum function ~\cite[Theorem 5.7]{Ten} $\sigma(j) \ll j \log \log j$,  we have 
\[
g(j) = \sum_{h | j} h = \sigma(j) \ll j\log\log j,\quad 1 \leq j < b.
\]
By the bound on the divisor function ~\cite[Theorem 5.4]{Ten}
\[
d (j) \ll f(j) = \exp\left(\frac{(1 + \varepsilon) (\ln 2) \log j}{\log \log j } \right),
\]
 we have
\[
g(j) = \sum_{\substack{h|j\\1\leq h<b}}h \leq b d(j) \ll b f(j),\quad b \leq j < b H.
\] 
Noting that $f$ is an increasing function, we have
\begin{align}\label{I}
I =& \sum_{1\le j<b H}
\min\left(2b,\left\|\frac{2\alpha j}a\right\|^{-1}\right)
\frac{g(j)}j \nonumber \\
=& \sum_{1 \le j<b}
\min\left(2b,\left\|\frac{2\alpha j}a\right\|^{-1}\right)
\frac{\sigma(j)}j
+ \sum_{b \le j<b H}
\min\left(2b,\left\|\frac{2\alpha j}a\right\|^{-1}\right)
\frac{g(j)}j \nonumber \\
\ll& (\log \log b) \sum_{1 \le j<b}
\min\left(2b,\left\|\frac{2\alpha j}a\right\|^{-1}\right)
+ b \sum_{b \le j<b H}
\min\left(2b,\left\|\frac{2\alpha j}a\right\|^{-1}\right)
\frac{f(j)}j \nonumber \\
\ll& (\log \log b) 
\sum_{1 \le j < b}
\min\left(2b,\left\|\frac{2\alpha j}a\right\|^{-1}\right)
+ b f(b H) \sum_{b \le j<b H}
\min\left(2b,\left\|\frac{2\alpha j}a\right\|^{-1}\right)
\frac{1}{j}. 
\end{align}
Hence,
we are left with the following two sums
\[I_1 = \sum_{1 \le j < b}
\min\left(2b,\left\|\frac{2\alpha j}a\right\|^{-1}\right)\]
and 
\[I_2 = \sum_{b \le j<b H}
\min\left(2b,\left\|\frac{2\alpha j}a\right\|^{-1}\right)
\frac{1}{j}.
\]

Let $q = \frac{a}{2 \alpha}$. When $\frac{j}{q}$ is close to an integer, the term $\left\|\frac{j}{q}\right\|^{-1}$ may be larger than $2b$. Precisely, for $1 \leq j < b H$, we have $\min\left(2b,\left\|\frac{j}{q}\right\|^{-1}\right) = 2b$ if $1 \leq j \leq \frac{q}{2b}$, or if there exists some integer $1 \leq k \leq \frac{b H}{q}$, such that $k - \frac{1}{2 b} \leq \frac{j}{q} \leq k + \frac{1}{2 b}$. 
For other $j$ on the interval $[1, b H)$, i.e., for $\frac{q}{2b} < j \leq \frac{q}{2}$ or $\left(k + \frac{1}{2b}\right) q < j < \left(k + 1 - \frac{1}{2b}\right)q$ for some integer $1 \leq k \leq \frac{b H}{q}$, we have $\min\left(2b,\left\|\frac{j}{q}\right\|^{-1}\right) = \left\|\frac{j}{q}\right\|^{-1}$. Moreover, we have $\left\|\frac{j}{q}\right\|^{-1} = \frac{q}{j - k q}$ when $\left(k + \frac{1}{2b}\right) q < j \leq \left(k + \frac{1}{2}\right)q$ for any integer $1 \leq k \leq \frac{b H}{q}$, since $\frac{j}{q}$ is closer to the integer $k$ than to $k + 1$. Similarly, we have $\left\|\frac{j}{q}\right\|^{-1} = \frac{q}{k q - j}$ when $\left(k - \frac{1}{2}\right)q < j < \left(k - \frac{1}{2b}\right)q$ for any integer $1 \leq k \leq \frac{b H}{q}$, since $\frac{j}{q}$ is closer to the integer $k$ than to $k - 1$. 

Note that when $b^2 \leq a$, Theorem~\ref{t1} follows from the trivial bound $E_\mathcal{P}(a, b)\ll b \ll \sqrt{a}$. For the rest of the proof, we will assume $b > \sqrt{a}$. Thus 
\begin{align}\label{I1}
I_1 
= & \sum_{1 \le j < b}
\min\left(2b,\left\|\frac{2\alpha j}a\right\|^{-1}\right) \nonumber \\
\ll &
\sum_{1 \leq j \leq \frac{q}{2b}} b
+ \sum_{\frac{q}{2 b} < j \leq \frac{q}{2}} \frac{q}{j}
+ 
\sum_{1 \leq k \le \frac{b}{q}}
\left(
\sum_{\left(k - \frac{1}{2b} \right)q \leq j \leq \left(k + \frac{1}{2 b} \right)q} 
b\right.\nonumber\\
&\phantom{x}\left.+\sum_{\left(k+\frac{1}{2b}\right)q < j \leq \left(k+\frac{1}{2}\right)q}
\frac{q}{j - k q} 
+\sum_{\left(k-\frac{1}{2}\right)q < j < \left( k - \frac{1}{2b} \right)q}
\frac{q}{k q - j} 
\right) \nonumber \\
\ll& b \left( \frac{q}{2 b} + 1 \right) + q \log q
+ b \left( \frac{q}{b} + 1\right) \left( \frac{b}{q} + 1 \right)
+ \sum_{1 \leq k \leq \frac{b}{q}}
2 q \sum_{\frac{q}{2b} < \ell \leq \frac{q}{2}}
\frac{1}{\ell} \nonumber \\
\ll& a + b + a \log a
+ \frac{a^2 + 2 a b + b^2}{a}
+ a \left( \frac{b}{a} + 1\right) \log a \nonumber \\ 
\ll& \frac{b^2}{a} + (a + b) \log a \nonumber \\
\ll& \frac{a^2 + b^2}{a} \log a.
\end{align}
We may treat $I_2$ in a similar fashion.
\begin{align}\label{I2}
I_2 
=& \sum_{b \le j< b H}
\min\left(2b,\left\|\frac{2\alpha j}a\right\|^{-1}\right)
\frac{1}{j} \nonumber \\
\ll& \sum_{\substack{1 \leq j \leq \frac{q}{2b} \\ j \geq b}} \frac{b}{j}+
\sum_{\substack{\frac{q}{2 b} < j \leq \frac{q}{2} \\ j \geq b}} \frac{q}{j^2}\nonumber
+ \sum_{1 \leq k \le \frac{b H} q}
\left(
\sum_{\substack{\left(k - \frac{1}{2b} \right)q \leq j \leq \left(k + \frac{1}{2 b} \right)q \\ j \geq b}} 
\frac{b}{j}\right.\\
&\phantom{x}+ \left.\sum_{\substack{\left(k+\frac{1}{2b}\right)q < j \leq \left(k+\frac{1}{2}\right)q \\ j \geq b}}
\frac{q}{j - k q}\frac{1}{j}
+ \sum_{\substack{\left(k-\frac{1}{2}\right)q < j < \left( k - \frac{1}{2b} \right)q \\ j \geq b}}
\frac{q}{k q - j}\frac{1}{j}
\right) \nonumber \\
\ll& 
\frac{q}{b}+1
+ \frac{b}{q} \left( \frac{q}{b} + 1\right) \log \left( \frac{b H}{q} + 2\right)\nonumber\\
&\phantom{x}+ \sum_{1 \leq k \leq \frac{bH}{q}}
\left(
\sum_{\substack{\left(k+\frac{1}{2b}\right)q < j \leq \left(k+\frac{1}{2} \right)q \\ j \geq b}}
\frac{q}{j- k q}\frac{1}{j}
+\sum_{\substack{\left(k-\frac{1}{2}\right)q < j < \left( k - \frac{1}{2b} \right) q  \\ j \geq b}}
\frac{q}{k q - j}\frac{1}{j}
\right) \nonumber \\
\ll& 
\frac{q}{b}
+ \left( \frac{b}{q} + 1\right) \log \left( \frac{b H}{q} + 2\right)\nonumber\\
&\phantom{x}+ \sum_{1 \leq k \le \frac{b H} q} \frac{1}{k - \frac{1}{2}} 
\left(
\sum_{\substack{\left(k+\frac{1}{2b}\right)q < j \leq \left(k+\frac{1}{2} \right)q \\ j \geq b}}
\frac{1}{j - k q}
+ \sum_{\substack{\left(k-\frac{1}{2}\right)q < j < \left( k - \frac{1}{2b} \right) q  \\ j \geq b}}
\frac{1}{k q - j}
\right) \nonumber \\ 
\ll& 
\left( \frac{b}q + \frac{q}b\right) \log \left( \frac{b H}{q} + 2\right)
+ \log\left(\frac{b H}{q} + 2\right) 
\sum_{\frac{q}{2b} < \ell \leq \frac{q}{2}}
\frac{1}{\ell} \nonumber \\ 
\ll& 
\left( \frac{b}{q} + \frac{q}b\right) \log \left( \frac{b H}{q} + 2\right)
+ \log\left(\frac{b H}{q} + 2\right) \log q \nonumber \\ 
\ll& \frac{a^2 + b^2}{a b} \log\left(\frac{b H}{a} + 2\right) \log a. 
\end{align}
We take $H = \left\lfloor \frac{2\sqrt{a}b }{a + b} \right\rfloor$. Since $a>1$ and $b>\sqrt{a}$, it is easily verified that $H\ge1$ and $\log bH\gg \log ab$. Then it follows from \eqref{I}, \eqref{I1}, and \eqref{I2} that
\begin{align}\label{Ibound}
I =& (\log \log b) I_1 + b f(b H) I_2 \nonumber \\
\ll& (\log \log b) \frac{a^2 + b^2}{a} \log a
+ \frac{a^2 + b^2}{a} f(b H) \log\left(\frac{b H}{a} + 2\right) \log a \nonumber \\
\ll& \left(a+\frac{b^2}{a}\right) 
\exp\left( \frac{\left(\ln 2 + \varepsilon \right) \log \Delta}{\log \log \Delta} \right), 
\end{align}
where $\Delta =\max\left( \frac{b^2 a^{\frac{1}{2}}}{a + b},3\right)$. 

Since $\frac{b}{H} \ll a^{\frac{1}{2}} + b a^{-\frac{1}{2}}$ and $b^{\frac12} \ll a^{\frac12}+ba^{-\frac12}$, by \eqref{e2.2} and \eqref{Ibound} we have
$$
E_{\mathcal{P}} (a, b)
\ll \left(a^{\frac{1}{2}} + b a^{- \frac{1}{2}}\right) \exp\left( \frac{\left(\frac{1}{2}\ln 2 + \varepsilon \right) \log \Delta}{\log \log \Delta } \right). 
$$
This completes the proof. 

\section{The proof of Theorem \ref{t2}}

The similarity between Theorem \ref{t1} and \ref{t2} is more than superficial. Indeed, seasoned workers in analytic number theory should have no difficulty realizing that both questions boil down to estimating the same type of exponential sums. Therefore, the proof Theorem \ref{t2} is almost identical to that of Theorem \ref{t1}, except that we replace Vaaler's lemma by the Erd\H{o}s-Turan inequality, which we introduce below for completeness. 

 Let $\{u_n\}_{1 \leq n \leq N}$ be a sequence of $N$ real numbers and $Z(N;\xi,\eta)$ count the number of $u_n$ such that $u_n\in(\xi,\eta)\bmod{1}$, where $\xi<\eta<\xi+1$. The discrepancy function 
\begin{equation*}
D(N;\xi,\eta)=Z(N;\xi,\eta)-(\eta-\xi)N
\end{equation*}
has been studied extensively, and we refer the readers to \cite[Chapter 1]{Mo} for more background. 
\begin{lemma}[{\cite[Corollary 1.1]{Mo}}]\label{Lemma2}
For any positive integer $H$, we have
\[
D(N; \xi, \eta) \leq \frac{N}{H + 1}
+ 3\sum_{1 \leq h \leq H}
\frac{1}{h}
\left|
\sum_{1 \leq n \leq N}
e(h u_n)
\right|. 
\]
\end{lemma}
Recall that  
\[
A_\mathcal{P}(a, b, \delta) 
= \sum_{\substack{1 \leq x \leq b \\ \left\| \frac{\alpha x^2}a+\beta x+\gamma a \right\| < \delta }} 1.
\]
To prove Theorem \ref{t2}, we apply Lemma \ref{Lemma2} to the sequence $\left\{\frac{\alpha x^2}{a} + \beta x + \gamma a \right\}_{1 \leq x \leq b}$ with  $\xi = - \delta$, $\eta = \delta$, and obtain
\begin{align*}
    A_\mathcal{P}(a, b, \delta) -2\delta \lfloor b\rfloor
   \ll\frac{b}{H+1}+3\sum_{1\leq h \leq H}\frac1h\left|\sum_{1\leq x\le b}e\left( h \left( \frac{\alpha x^2}{a} + \beta x + \gamma a \right) \right)\right|.
\end{align*}
Now we proceed in the same way as the treatment of the right hand side of \eqref{e2.1} in the previous section, and therefore complete the proof of Theorem \ref{t2}. 

\proof[Acknowledgement]
The authors are grateful to the anonymous referee for helpful suggestions and comments. 

\bibliographystyle{plain}
 \bibliography{biblio}

\end{document}